\documentclass[fleqn]{mat01}
\usepackage{times,mathtimy,amssymb,latexsym}
\begin{document}

\setcounter{page}{493} \firstpage{493}

%%%%%%%%%%%new commands %%%%%%%%%%%%%%%%%
\newcommand{\IA}{{\Bbb A}} \newcommand{\IB}{{\Bbb B}}
\newcommand{\IC}{{\Bbb C}} \newcommand{\ID}{{\Bbb D}}
\newcommand{\IE}{{\Bbb E}} \newcommand{\IF}{{\Bbb F}}
\newcommand{\IG}{{\Bbb G}} \newcommand{\IH}{{\Bbb H}}
\newcommand{\II}{{\Bbb I}} \newcommand{\IJ}{{\Bbb J}}
\newcommand{\IK}{{\Bbb K}} \newcommand{\IL}{{\Bbb L}}
\newcommand{\IM}{{\Bbb M}} \newcommand{\IN}{{\Bbb N}}
\newcommand{\IO}{{\Bbb O}} \newcommand{\IP}{{\Bbb P}}
\newcommand{\IQ}{{\Bbb Q}} \newcommand{\IR}{{\Bbb R}}
\newcommand{\IS}{{\Bbb S}} \newcommand{\IT}{{\Bbb T}}
\newcommand{\IU}{{\Bbb U}} \newcommand{\IV}{{\Bbb V}}
\newcommand{\IW}{{\Bbb W}} \newcommand{\IX}{{\Bbb X}}
\newcommand{\IY}{{\Bbb Y}} \newcommand{\IZ}{{\Bbb Z}}

\newtheorem{theore}{Theorem}
\renewcommand\thetheore{\arabic{section}.\arabic{theore}}
\newtheorem{theor}[theore]{\bf Theorem}
\newtheorem{lem}[theore]{\it Lemma}
\newtheorem{propo}[theore]{\rm PROPOSITION}

\title{A functional central limit theorem for a class of urn models}

\markboth{Gopal K Basak and Amites Dasgupta}{A functional central
limit theorem for a class of urn models}

\author{GOPAL K BASAK$^{1,2}$ and AMITES DASGUPTA$^{2}$}

\address{$^{1}$Department of Mathematics, University of Bristol, University
Walk, Bristol, BS8~1TW, UK\\
\noindent $^{2}$Stat-Math Unit, Indian Statistical
Institute, 203 B.T.~Road, Kolkata~700~108, India\\
\noindent E-mail: magkb@bris.ac.uk and gkb@isical.ac.in;
amites@isical.ac.in}

\volume{115}

\mon{November}

\parts{4}

\pubyear{2005}

\Date{MS received 31 March 2005}

\begin{abstract}
We construct an independent increments Gaussian process associated
to a class of multicolor urn models. The construction uses random
variables from the urn model which are different from the random
variables for which central limit theorems are available in the
two color case.
\end{abstract}


\maketitle

\section{Introduction}

Consider a four-color urn model in which the replacement matrix is
actually a stochastic matrix $\mathbf{R}$ as in
ref.~\cite{gouet2}. That is, we start with one ball of any color,
which is the \hbox{0-th} trial. Let $\mathbf{W}_n$ denote the
column vector of the number of balls of the four colors up to the
$n$-th trial, where the components of $\mathbf{W}_n$ are
nonnegative real numbers.  Then a color is observed by random
sampling from a multinomial distribution with probabilities
$(1/(n+1)) \mathbf{W}_n$. Depending on the color that is observed,
the corresponding row of $\mathbf{R}$ is added to
$\mathbf{W}_n^\prime$ and this gives $\mathbf{W}_{n+1}^\prime$. A
special case of the main theorem of Gouet \cite{gouet2} is that if
the stochastic matrix $\mathbf{R}$ is irreducible, then $(1/(n+1))
\mathbf{W}_n^\prime$  converges almost surely (a.s.) to the
stationary distribution $\mathbf{\pi}$ of the irreducible
stochastic matrix $\mathbf{R}$ (it should be carefully noted that
the multicolor urn model is vastly different from the Markov chain
evolving according to the transition matrix equal to the
stochastic matrix $\mathbf{R}$, also notice that $\pi$ is a row
vector). Suppose the nonprincipal eigenvalues of $\mathbf{R}$
satisfy $\lambda_1 < 1/2, \lambda_2 = 1/2, \lambda_3
> 1/2$ respectively, which are assumed to be real (and hence lie
in $(-1, 1)$), and $\xi_1, \xi_2, \xi_3$ be the corresponding
eigenvectors. Using $\mathbf{\pi} \xi_i = \mathbf{\pi} \mathbf{R}
\xi_i = \lambda_i \mathbf{\pi} \xi_i$ it is seen that $(1/(n+1))
\mathbf{W}_n^\prime \xi_i \rightarrow 0$.

Central and functional central limit theorems for
$\mathbf{W}_n^\prime \xi_i$ have been the subject of several
papers in the literature \cite{freedman,gouet,smythe} especially
for two-color models and also some multicolor models. The norming
in the central limit theorems in the two color urn models depends
on the nonprincipal eigenvalue as follows: for $\lambda < 1/2$ the
rate is $\sqrt{n}$, for $\lambda = 1/2$ the rate is $\sqrt{n \log
n}$ and the limits are normal in these two cases. However for
$\lambda > 1/2$ the  rate is $\Pi_0^{n - 1} (1 + (\lambda/(j+1))$
and in this case the limit exists almost\break surely.

Functional central limit theorems (FCLT) for a class of two-color
urn models have been considered by Gouet \cite{gouet}. These
FCLT's of Gouet \cite{gouet} use the same norming,  as stated in
the previous paragraph, under which central limit theorems have
been proved in \cite{freedman} and \cite{gouet}.
Ref.~\cite{janson} contains a survey of the literature on such
FCLT's. In this article we prove a different FCLT that uses random
variables with the norming $\Pi_0^{n - 1} (1 + (\lambda/(j+1))$
irrespective of whether $\lambda$ is less than 1/2, equal to 1/2
or greater than 1/2. This is the main result of the paper. For the
sake of convenience we restrict ourselves to real eigenvalues
only. We state the result for the above four-color model but it
can be seen from the proof that it can be extended to urn models
with any number of\break colors.

The article is organized as follows. In \S2 we develop the
notation, state the main result and give its proof. Some of the
calculations have been done separately in\break \S3.

\section{Main result}

We write
\begin{equation}
Z_{i,n} =  \frac{\mathbf{W}_n^\prime \xi_i}{\Pi_0^{n - 1} \left(1
+ \frac{\lambda_i}{j+1}\right)},\label{fera}
\end{equation}
where $\xi_i$ is the eigenvector corresponding to the eigenvalue
$\lambda_i$. From the description of the urn model we have
$\mathbf{W}_{n+1}^\prime \xi_i = \mathbf{W}_n^\prime \xi_i +
\chi_{n+1}^\prime \mathbf{R} \xi_i = \mathbf{W}_n^\prime \xi_i +
\lambda_i \chi_{n+1}^\prime \xi_i$, where $\chi_{n+1}$ is the
column vector consisting of the indicator functions of balls of
the four colors respectively. We also have
\begin{equation}
E \{ \chi_{n+1}^\prime \xi_i |{\cal F}_n \} = \frac{1}{n+1}
\mathbf{W}_n^\prime \xi_i,\label{hir}
\end{equation}
where ${\cal F}_n$ is the $\sigma$-field of observations up to the
$n$-th trial. From this it follows that $Z_{i,n}$ is a martingale.
From \S3, it follows that $Z_{3,n}$ is $L^2$-bounded so that it
converges almost surely. However in the two color case, for
$\lambda < 1/2, \mathbf{W}_n^\prime \xi/\sqrt{n}$ and for $\lambda
= 1/2, \mathbf{W}_n^\prime \xi/ \sqrt{n \log n}$ converge to
normal distributions and the FCLT's proved in \cite{gouet} use
such normalizations. Thus the question of using the same norming
$\Pi_0^{n - 1} (1 + (\lambda/(j+1)) \sim n^{\lambda}$, to get an
FCLT irrespective of $\lambda < 1/2, \lambda = 1/2$ or $\lambda >
1/2$, is of interest. Our main result, Proposition~2.1, is a step
in this direction using the tails of the sequence $(Z_{1,n}, Z_{2,
n}, Z_{3, n})$ whereas the FCLT's in the literature are based on
partial sums starting from the\break beginning.

\begin{propo}$\left.\right.$\vspace{.5pc}

\noindent The sequence of processes $\mathbf{G}_n(t) =
(G_{1,n}(t), G_{2,n}(t), G_{3,n}(t))$ where
\begin{equation*}
G_{i, n}(t) = \sum_{m = n}^{[ne^t]} m^{\lambda_i - 1/2} (Z_{i,
m+1} - Z_{i, m}),\quad i = 1, 2, 3, \ t \geq 0,
\end{equation*}
converges to an independent increments Gaussian process
$\mathbf{G}(t)$ with covariance function $c_{i,j}(t) = t \lambda_i
\lambda_j \langle \xi_i \xi_j, \pi^\prime \rangle, i, j = 1, 2,
3${\rm ,} where the vector of the coordinate-wise product of the
components of the two vectors $\xi_i$ and $\xi_j$ is denoted by
$\xi_i \xi_j$ and the Euclidean inner product of the two vectors
is denoted by $\langle . , .\rangle $.
\end{propo}
Note that the process $\mathbf{G}$ can be viewed as a
multidimensional Wiener process with covariance function $c_{i,
j}(\cdot)$.

\begin{proof}
From eq.~(\ref{fera}) we have the following expansion:
\begin{align}
Z_{i, m+1} - Z_{i,m} &\sim  - \frac{\lambda_i}{m} Z_{i,m} +
\lambda_i \frac{\chi_{m+1}^\prime \xi_i}{\Pi_{k=0}^{m} \left(1 +
\frac{\lambda_i}{k+1}\right)}\nonumber\\[.2pc]
&\sim - \frac{\lambda_i}{m} Z_{i,m} +  \lambda_i
\frac{\chi_{m+1}^\prime \xi_i}{m^{\lambda_i}}.\label{hi}
\end{align}
Since the components of $\mathbf{G}_n(t)$ are martingales, an
independent increments Gaussian process as a limiting process is
expected. In particular we follow Theorem~1.4, p.~339 of
\cite{ek}, by which it is enough to show that the joint
characteristics of the martingales converge to a joint covariance
function. Note that from (\ref{hi}) $m^{\lambda_i -
\frac{1}{2}}(Z_{i,m+1} - Z_{i, m}) = O(1/\sqrt{m})$, as
$\mathbf{W}_m^\prime \xi_i/m$ and $\chi_{m+1}^\prime \xi_i$ are
bounded. This takes care of continuity of the paths and cross
quadratic variations which is condition (b) of that theorem. Thus
it remains to show that the cross quadratic variations converge to
$c_{i,j}(t)$. We first do this for $i = 1, j = 2$. From (\ref{hi})
we have
\begin{align}
m^{\lambda_1 - 1/2}(Z_{1,m+1} - Z_{1,m}) &\sim - \lambda_1
m^{\lambda_1 - 1/2} \frac{\mathbf{W}_m^\prime \xi_1}{m^{\lambda_1
+ 1}} + \lambda_1 m^{\lambda_1 - 1/2} \frac{(\chi_{m+1}^\prime
\xi_1)}{m^{\lambda_1}},\nonumber\\[.2pc]
m^{\lambda_2 - 1/2}(Z_{2,m+1} - Z_{2,m}) &\sim - \lambda_2
m^{\lambda_2 - 1/2} \frac{\mathbf{W}_m^\prime \xi_2}{m^{\lambda_2
+ 1}} + \lambda_2 m^{\lambda_2 - 1/2} \frac{(\chi_{m+1}^\prime
\xi_2)}{m^{\lambda_2}}.\label{fi}
\end{align}
We want to show that in computing the cross quadratic variation,
which is the limit of
\begin{equation*}
\sum_n^{[n e^t]} E\{ m^{\lambda_1 - 1/2} m^{\lambda_2 - 1/2}
(Z_{1, m+1} - Z_{1, m})(Z_{2, m+1} - Z_{2, m})|{\cal F}_m  \},
\end{equation*}
only the second term from the right-hand side of each of
eqs~(\ref{fi}) contributes. Since $\chi_{n+1}$ consists of
indicator functions, which implies that
\begin{equation*}
\left(\sum_k \xi_{1,k} \chi_{n+1, k}\right)\left(\sum_l \xi_{2,l}
\chi_{n+1, l}\right) = \sum_k \xi_{1,k} \xi_{2, k} \chi_{n+1, k},
\end{equation*}
this contribution is the limit of
\begin{equation*}
\lambda_1 \lambda_2 \sum_n^{[ne^t]} \frac{1}{m} \left\langle \xi_1
\xi_2, \frac{\mathbf{W}_m}{m+1} \right\rangle,
\end{equation*}
which is $t \lambda_1 \lambda_2 \langle \xi_1 \xi_2, \pi^\prime
\rangle$, since from \cite{gouet2} we know
${\mathbf{W}^\prime_m}/{(m+1)} \rightarrow \pi$ a.s. Also notice
that this part of the argument does not depend on whether
$\lambda_1$ or $\lambda_2$ are less than or equal to 1/2.

To see why the contribution to the cross quadratic variation from
the first terms of (\ref{fi}) goes to $0$, by Cauchy--Schwarz
inequality it is enough to show that the sum of squares over $n$
to $[n e^t]$ of the first terms in each line of (\ref{fi}) goes to
$0$. This part of the argument will depend on the value of
$\lambda_i$. Note the following which have been proved in \S3:
\begin{align}
\hbox{For} \ \lambda_1 &< 1/2,\quad \frac{\mathbf{W}_m^\prime
\xi_1}{\sqrt{m}} \mbox{ is } L^2\hbox{-bounded},\\[.3pc]
\hbox{For} \ \lambda_2 &= 1/2,\quad \frac{\mathbf{W}_m^\prime
\xi_2}{\sqrt{m\log m}} \mbox{ is } L^2\hbox{-bounded}.
\end{align}
Consider the case $\lambda_1 < 1/2$. We need to show
\begin{equation*}
\sum_n^{[n e^t]} \frac{( \mathbf{W}_m^\prime \xi_1)^2}{m^3}
\rightarrow  0\ \hbox{a.s}.
\end{equation*}
We know that for $\lambda_1 < 1/2$, $\mathbf{W}_m^\prime \xi_1
/\sqrt{m}$ is $L^2$-bounded, so that
\begin{equation}
E \sum_n^{[n e^t]} \frac{( \mathbf{W}_m^\prime \xi_1)^2}{m^3} \leq
E \sum_n^\infty  \frac{( \mathbf{W}_m^\prime \xi_1)^2}{m^3} \leq
\hbox{const.} \sum_n^\infty \frac{1}{m^2} \rightarrow
0.\label{fii}
\end{equation}
Since the sum inside the expectation in the middle is decreasing
in $n$, it converges to $0$ a.s. For $\lambda_2 = 1/2$,
$\mathbf{W}_m^\prime \xi_2 /\sqrt{m\log m}$ is $L^2$-bounded, and
one can proceed similarly. Thus we have proved that $c_{1, 2}(t) =
t \lambda_1 \lambda_2 \langle \xi_1 \xi_2, \pi^\prime \rangle$.
Similarly $c_{i, j}(t), i, j = 1, 2$, can be computed as given in
Proposition~2.1.

Now consider as to what will happen if we were computing say
$c_{1,3}(t)$. For $\lambda_3 > 1/2$, the expansion (\ref{fi}) is
similar, and in the cross quadratic variation the contribution of
the second term from the right-hand side of
\begin{equation*}
m^{\lambda_3 - 1/2}(Z_{3,m+1} - Z_{3,m}) \sim - \lambda_3
m^{\lambda_3 - 1/2} \frac{\mathbf{W}_m^\prime \xi_3}{m^{\lambda_3
+ 1}} + \lambda_3 m^{\lambda_3 - 1/2} \frac{(\chi_{m+1}^\prime
\xi_3)}{m^{\lambda_3}}
\end{equation*}
is similar to what we had before. For $\lambda_3 > 1/2$,
$\mathbf{W}_m^\prime \xi_3 / \Pi_0^{m - 1} (1 + (\lambda_3/(j+1))$
is a martingale and from Appendix~3.3,
\begin{equation}
\frac{\mathbf{W}_m^\prime \xi_3}{m^{\lambda_3}} \hbox{ is } L^2
\hbox{-bounded}.
\end{equation}
So $\mathbf{W}_m^\prime \xi_3 / m^{\lambda_3}$ converges almost
surely. This implies that the contribution of the first term
\begin{equation*}
\sum_n^{[ne^t]} \frac{(\mathbf{W}_m^\prime \xi_3)^2}
{m^{2\lambda_3}} \frac{1}{m^{3 - 2\lambda_3}} \rightarrow 0 \hbox{
a.s. }
\end{equation*}
since $2 \lambda_3 < 2$. Thus $c_{i, j}(t), i = 1,2,3, j = 3$, can
be computed as given in the statement of Proposition~2.1. This
completes the proof.\hfill $\Box$
\end{proof}

\section{Appendix}

Suppose real eigenvalues satisfy $\lambda_1 < 1/2, \lambda_2 =
1/2, \lambda_3 > 1/2$ and $\xi_1, \xi_2, \xi_3$ be the
corresponding eigenvectors. In this section we prove that $X_n$,
$Y_n$ and $Z_n$ are $L^2$-bounded where
\begin{equation}
X_n = \frac{\mathbf{W}_n^\prime \xi_1}{\sqrt{n}},\quad Y_n =
\frac{\mathbf{W}_n^\prime \xi_2}{\sqrt{n \log n}},\quad Z_n =
\frac{\mathbf{W}_n^\prime \xi_3}{\Pi_0^{n - 1} \left(1 +
\frac{\lambda_3}{j+1}\right)},
\end{equation}
a fact which has been used in the proof of Proposition~2.1. For
$X_n$ and $Y_n$ verification of $L^2$-boundedness is through
Lemma~2.1 of \cite{kersting}. This is done on a case by case basis
depending on $\lambda_1$ and $\lambda_2$ in the next two
subsections. For the reader's convenience we state Kersting's
lemma from \cite{kersting} here:

\setcounter{section}{2} \setcounter{theore}{0}
\begin{lem}\hskip -.4pc{\rm \cite{kersting}.} \ \ Let $\alpha_n,
\beta_n\ (n \geq 1)$ be nonnegative numbers such that $\alpha_n
\rightarrow 0${\rm ,} $\sum_{n = 1}^\infty \alpha_n = \infty$, and
for large $n${\rm ,}
\begin{equation*}
\beta_{n+1} \leq \beta_n ( 1 - c \alpha_n) + d \alpha_n
\end{equation*}
with $c, d > 0$. Then $\limsup_{n \rightarrow \infty} \beta_n \leq
d/c$.
\end{lem}
\setcounter{section}{3}

\subsection{\it $L^2$-boundedness of $X_n$}

Using $\mathbf{W}_{n+1}^\prime \xi_1 = \mathbf{W}_n^\prime \xi_1 +
\lambda_1 \chi_{n+1}^\prime \xi_1$ and the definition of $X_n$, we
get
\begin{equation*}
X_{n+1} =  X_n \sqrt{\frac{n}{n+1}} + \lambda_1
\frac{\chi_{n+1}^\prime \xi_1}{\sqrt{n+1}}.
\end{equation*}
Taking conditional expectation and using (\ref{hir}) we get
\begin{align*}
E \{ X_{n+1}^2|{\cal F}_n \} = X_n^2 \left(1 - \frac{1}{n+1}
\right) \left(1 + \frac{2\lambda_1}{n+1} \right) +
\frac{\lambda_1^2}{n+1} \left\langle \frac{\mathbf{W}_n}{n+1},
\xi_1^2\right\rangle,
\end{align*}
from which taking further expectation we get
\begin{equation*}
E X_{n+1}^2 = E X_n^2 \left(1 - \frac{1}{n+1}\right) \left(1 +
\frac{2\lambda_1}{n+1}\right) + \frac{\lambda_1^2}{n+1}
\left\langle E \frac{\mathbf{W}_n}{n+1}, \xi_1^2 \right\rangle.
\end{equation*}
The last vector $E ({\mathbf{W}_n}/{(n+1)})$ consists of bounded
components. Thus if $\lambda_1 < 0$, then
\begin{equation*}
E X_{n+1}^2 \leq E X_n^2 \left(1 - \frac{1}{n+1}\right) +
\frac{{\rm const}}{n+1},
\end{equation*}
and Kersting's lemma applies. If $\lambda_1
> 0$ then we still have  $\lambda_1 < 1/2$ i.e. $2\lambda_1 < 1$.
In this case
\begin{align*}
\left(1 - \frac{1}{n+1}\right) \left(1 + \frac{2\lambda_1}{n+1}
\right) &\leq 1 +  \frac{2\lambda_1}{n+1} - \frac{1}{n+1}\\[.2pc]
&= 1 - \frac{1 - 2\lambda_1}{n+1},
\end{align*}
i.e. Kersting's lemma applies.

\subsection{\it $L^2$-boundedness of $Y_n$}

Using $\mathbf{W}_{n+1}^\prime \xi_2 = \mathbf{W}_n^\prime \xi_2 +
\lambda_2 \chi_{n+1}^\prime \xi_2$ and the definition of $Y_n$, we
get
\begin{equation*}
Y_{n+1} = Y_n \sqrt{\frac{n \log n}{(n+1) \log (n+1)}} + \lambda_2
\frac{\chi_{n+1}^\prime \xi_2}{\sqrt{(n+1) \log (n+1)}}.
\end{equation*}
Taking conditional expectation we get (recall $\lambda_2 = 1/2$)
\begin{align*}
E \{ Y_{n+1}^2|{\cal F}_n \} &= Y_n^2 \frac{n \log n}{(n+1) \log
(n+1)} \left(1 + \frac{1}{n+1} \right)\\[.2pc]
&\quad\, + \frac{\lambda_2^2}{(n+1) \log (n+1)} \left\langle
\frac{\mathbf{W}_n}{n+1}, \xi_2^2 \right\rangle,
\end{align*}
from which taking further expectation we get
\begin{align*}
E Y_{n+1}^2 &=  E Y_n^2 \left( 1 - \frac{(n+1) \log (n+1) - n \log
n}{(n+1) \log (n+1)} \right) \left(1 + \frac{1}{n+1} \right)\\[.2pc]
&\quad\, + \frac{\lambda_2^2}{(n+1) \log (n+1)} \left\langle E
\frac{\mathbf{W}_n}{n+1}, \xi_2^2 \right\rangle.
\end{align*}
The vector $E ({\mathbf{W}_n}/{(n+1)})$ consists of bounded
components. Now we apply the second trick of the previous
subsection to apply Kersting's lemma. The following calculation
does the rest of the work.

We show that
\begin{equation*}
(n+1)\log (n+1) \left\{ \frac{(n+1)\log (n+1) - n\log n}{(n+1)\log
(n+1)} - \frac{1}{n+1} \right\} \rightarrow c > 0.
\end{equation*}
We approximate $\log (n+1)$ by $\log n + \frac{1}{n}$. This gives
\begin{align*}
&(n+1) \log (n+1) - n\log n \sim (n+1) \log n + (n+1)\frac{1}{n} -
n\log n\\[.2pc]
&\quad\, = \log n + 1 + \frac{1}{n}.
\end{align*}
Hence\vspace{-.2pc}
\begin{align*}
\left(\log n + 1 + \frac{1}{n} \right) - \log (n+1) \sim
\left(\log n + 1 + \frac{1}{n} \right) - \left(\log n +
\frac{1}{n}\right) \rightarrow 1.
\end{align*}

\vspace{.1pc}
\subsection{\it $L^2$-boundedness of $Z_n$}

The proof follows Lemma~3.1 of \cite{freedman}. We have earlier
proved the approximation
\begin{equation*}
Z_{n+1} - Z_n \sim - \frac{\lambda_3}{n} Z_n + \lambda_3
\frac{\chi_{n+1}^\prime \xi_3}{n^{\lambda_3}}.
\end{equation*}
Now with the martingale property of $Z_n$, $Z_{n+1} = Z_n +
(Z_{n+1} - Z_n)$, and by the above approximation we have
\begin{equation}
E (Z_{n+1}^2|{\cal F}_n) \sim Z_n^2 +
\left(\frac{\lambda_3^2}{n^2} Z_n^2 - 2 \frac{\lambda_3^2}{n^2}
Z_n^2 + \lambda_3^2 \frac{\left\langle \frac{\mathbf{W}_n}{n+1},
\xi_3^2 \right\rangle}{n^{2\lambda_3}} \right),\nonumber
\end{equation}
implying
\begin{equation}
E Z_{n+1}^2 \leq E Z_n^2 \left( 1 - \frac{\lambda_3^2}{n^2}
\right) + \lambda_3^2 \frac{{\rm const}}{n^{2\lambda_3}}.
\label{hii}
\end{equation}
Since $2 \lambda_3 > 1$, by iteration of (\ref{hii}) it follows
that $Z_n$ is $L^2$-bounded.

\end{document}